\documentclass[11pt]{article}
\usepackage[letterpaper,margin=1in]{geometry}
\usepackage[T1]{fontenc}
\usepackage{lmodern}
\usepackage{amsmath,amssymb,amsthm,mathtools}
\usepackage{booktabs,array}
\usepackage{graphicx}
\usepackage{xcolor}
\usepackage{enumitem}
\usepackage{microtype}
\usepackage{tikz}
\usetikzlibrary{arrows.meta,positioning,fit}
\usepackage{float}
\usepackage{caption}
\usepackage{authblk}

\usepackage[colorlinks=true,linkcolor=blue!45!black,citecolor=blue!45!black,urlcolor=blue!45!black]{hyperref}

\newcommand{\R}{\mathbb R}
\newcommand{\NN}{\mathrm{NN}}
\newcommand{\Tt}[1]{[\![#1]\!]}
\newcommand{\TTSVD}{\operatorname{TT-SVD}}
\newcommand{\norm}[1]{\left\lVert #1\right\rVert}

\numberwithin{equation}{section}
\theoremstyle{plain}
\newtheorem{proposition}{Proposition}[section]
\theoremstyle{definition}
\newtheorem{definition}[proposition]{Definition}
\theoremstyle{remark}

\newtheorem{assumption}[proposition]{Assumption}

\setlist[itemize]{leftmargin=1.5em,itemsep=2pt,topsep=4pt}
\setlist[enumerate]{leftmargin=1.7em,itemsep=2pt,topsep=4pt}
\DeclareCaptionType[fileext=loa]{algorithm}[Algorithm][List of Algorithms]
\hypersetup{pdftitle={TEN-NO: Tensor Neural Operator},pdfauthor={Jinchen Wang, Rizheng Huang, Vladimir Kazeev, Yuehaw Khoo, Mathias Oster}}

\title{TEN-NO: Tensor Neural Operator}
\author[1]{Jinchen Wang\thanks{These authors contributed equally to this work.}}
\author[2]{Rizheng Huang\protect\footnotemark[1]}
\author[3]{Vladimir Kazeev}
\author[4]{Yuehaw Khoo}
\author[5]{Mathias Oster}
\affil[1]{Department of Statistics and Data Science, University of California, Los Angeles, Los Angeles, CA, USA}
\affil[2]{Oden Institute for Computational Engineering and Sciences, The University of Texas at Austin, Austin, TX, USA}
\affil[3]{Faculty of Mathematics, University of Vienna, Vienna, Austria}
\affil[4]{Department of Statistics and Committee on Computational and Applied Mathematics, University of Chicago, Chicago, IL, USA}
\affil[5]{Institute for Geometry and Practical Mathematics, RWTH Aachen University, Aachen, Germany}
\date{}

\begin{document}
\maketitle
\begingroup
\renewcommand{\thefootnote}{}
\footnotetext[0]{%
  \raggedright
  Author e-mail addresses:
  \href{mailto:jwang86@g.ucla.edu}{jwang86@g.ucla.edu}~(J.~Wang),
  \href{mailto:rh42322@my.utexas.edu}{rh42322@my.utexas.edu}~(R.~Huang),
  \href{mailto:vladimir.kazeev@univie.ac.at}{vladimir.kazeev@univie.ac.at}~(V.~Kazeev),
  \href{mailto:ykhoo@uchicago.edu}{ykhoo@uchicago.edu}~(Y.~Khoo),
  \href{mailto:oster@igpm.rwth-aachen.de}{oster@igpm.rwth-aachen.de}~(M.~Oster).
}
\endgroup

\begin{abstract}
We propose the Tensor Neural Operator (TEN-NO), a neural network architecture for predicting high-dimensional states whose dense representations suffer from the curse of dimensionality. TEN-NO represents states compactly in low-rank tensor-train (TT) form and learns directly in this representation, drawing on tensor-network methods for time evolution. We demonstrate TEN-NO's accuracy and storage efficiency in high dimensions by predicting states at a prescribed time from their initial conditions and by approximating solutions to parametric equations given their parameters. In particular, for Hamilton--Jacobi--Bellman (HJB) equations, Glauber spin dynamics, and transverse-field Ising models, TEN-NO achieves relative test errors of few percent or less while the dense tensor representations would contain up to $2^{256}$ entries. These results suggest that tensor networks provide a scalable foundation for operator learning in high dimensions.

\end{abstract}


\section{Introduction}\label{sec:introduction}
Learning a map between function spaces from data is a task of growing importance in scientific computing. Once learned, such a map replaces a numerical solver in many tasks, such as uncertainty quantification, inverse problems, and control. Operator learning has seen significant recent progress, particularly due to neural-network parametrizations, including networks for parametric PDEs~\cite{khoo2021parametric}, DeepONet~\cite{lu2021deeponet}, the Fourier neural operator (FNO)~\cite{li2021fno}, and neural operators more generally~\cite{kovachki2023neural}. All of these approaches represent the input and output functions by their values on a grid or at a set of points. This representation is effective for PDEs in two or three dimensions. In high dimensions, however, the size of a single input--output pair grows exponentially with the dimension $d$. For example, a function of $d$ variables with $32$ grid points per coordinate has $32^d$ values, and a probability distribution over $d$ spins has $2^d$ entries. We propose a neural-network architecture that learns maps between high-dimensional functions without storing the full tensors. 

More precisely, we consider high-dimensional states
\[
  \rho:[n_1]\times\cdots\times[n_d]\longrightarrow\R,
\]
such as probability distributions over $d$ spins, quantum wave functions, and functions of $d$ variables discretized on a tensor-product grid. In the \emph{parameter-to-state} setting, the target state $\rho_{\mathrm{tar}}$ is determined by parameters $\omega$, such as coefficients or external fields, through an equation

\begin{equation}
  \mathcal L_\omega(\rho_{\mathrm{tar}})=0.
  \label{eq:intro-param}
\end{equation}
For each sampled parameter vector $\omega^{(m)}$, a numerical solver computes the corresponding target state $\rho_{\mathrm{tar}}^{(m)}$. From the resulting $N$ pairs
\[
  \bigl(\omega^{(1)},\rho_{\mathrm{tar}}^{(1)}\bigr),\ \ldots,\ \bigl(\omega^{(N)},\rho_{\mathrm{tar}}^{(N)}\bigr),
  \qquad \mathcal L_{\omega^{(m)}}\bigl(\rho_{\mathrm{tar}}^{(m)}\bigr)=0,
\]
we learn the map $\omega\mapsto\rho_{\mathrm{tar}}$. In the \emph{evolution} setting, the state satisfies
\begin{equation}
  \partial_t\rho=\mathcal L_\omega(\rho),\qquad \rho(0)=\rho_0.
  \label{eq:intro-evol}
\end{equation}
This form occurs in master equations, such as the one for Glauber dynamics, and in HJB equations. We learn from initial states and their targets $\rho_{\mathrm{tar}}=\rho_T$, where $\rho_T:=\rho(T)$ at a fixed time $T$. Both settings use the prediction map $(\rho_0,\omega)\mapsto\rho_{\mathrm{tar}}$. Here $\rho_0$ is the initial state for evolution and a fixed reference state for a parameter-to-state map.

To learn these maps without ever forming the full tensors, we assume that the states admit low-rank TT representations~\cite{oseledets2011,kolda2009tensor,grasedyck2013survey}. Many operator-learning problems satisfy this assumption. For example, low-rank TTs accurately approximate probability distributions governed by master equations~\cite{johnson2010stochastic,temme2010stochastic,kazeev2014cme}, ground states of one-dimensional quantum Hamiltonians~\cite{White1993DMRG,schollwock2011}, and solutions of high-dimensional HJB equations~\cite{dolgov2021hjb,oster2022hjb}. In the TT format, a state with $d$ modes of size at most $n$ and ranks at most $r$ is stored in $O(dnr^2)$ numbers instead of $n^d$. The operations we need, such as addition, inner products, and rank truncation, act directly on the TT cores at a cost linear in $d$~\cite{oseledets2011}. We therefore propose TEN-NO, a neural network whose input, intermediate states, and output are all in TT form. Our contributions are as follows.
\begin{enumerate}
\item We develop a unified framework for operator learning in high dimensions that covers both fixed-time evolution maps and parameter-to-state maps. Since every state is kept in TT form and each layer updates the TT cores locally, the storage of the data and the costs of computing the updates and evaluating the loss all grow linearly in $d$.
\item Numerical experiments show that TEN-NO predicts terminal states of a $20$-dimensional HJB equation and of Glauber dynamics with $256$ spins, as well as ground states of the transverse-field Ising (TFI) model with $48$ sites, with test errors of a few percent or less. On a three-dimensional HJB problem, where FNO~\cite{li2021fno} can still be applied, our method reaches a test error comparable to that of FNO while providing a $17.7\times$ training speedup. It also compresses the inputs and reference solutions by a factor of $23.3$ relative to FNO's dense representation.
\end{enumerate}

\subsection{Related work}\label{sec:related}

\paragraph{Operator learning.}
Neural networks were used in early work to predict quantities of interest from the coefficients of parametric PDEs~\cite{khoo2021parametric}. DeepONet~\cite{lu2021deeponet}, FNO~\cite{li2021fno}, and neural operators more generally~\cite{kovachki2023neural} learn maps between function spaces. These architectures act on grid values or point evaluations. On a $d$-dimensional tensor-product grid with $n$ points per coordinate, a dense-grid FNO stores and processes $n^d$ values per channel. DeepONet evaluates its output at individual query points, with its branch network receiving samples of the input function at sensor locations. Recovering the output on the full grid requires $n^d$ point evaluations. Numerical structure has been used to reduce the number of parameters. Examples include SwitchNet~\cite{khoo2019switch}, multiscale networks based on hierarchical matrices~\cite{fan2019mnnh}, and the multi-grid tensorized FNO~\cite{kossaifi2023tfno}. The tensorized FNO factorizes the Fourier weights and operates on smaller grid patches while still storing the states as grid values. Tensor-train operator inference (TT-OpInf)~\cite{danis2025ttopinf} does keep the states in TT form. It fits polynomial dynamics to TT snapshots of trajectories and to estimated time derivatives, and integrates the fitted dynamics to make predictions. In contrast, we learn the map from input--output pairs alone, without trajectories or time derivatives. Our network predicts the full state in one forward pass, and every state along the way stays in TT form.

\paragraph{Neural-network solvers for high-dimensional problems.}
Neural networks are also used to solve individual high-dimensional problems. For PDEs, examples include the deep BSDE method~\cite{han2018bsde}, the deep Ritz method~\cite{e2018ritz}, and physics-informed neural networks~\cite{raissi2019pinn}. Neural-network quantum states~\cite{carleo2017nqs} are used to compute many-body ground states. In these methods, a network represents the solution of one instance of the equation and is trained by minimizing a residual or a variational objective. A new coefficient or initial condition requires a new training run, and the solution is accessed through point evaluations or samples. We instead learn the solution map over a family of inputs and return the full state as a TT.

\paragraph{Tensor-network solvers.}
In quantum physics, the TT format is known as the matrix product state (MPS)~\cite{schollwock2011}. The density matrix renormalization group (DMRG)~\cite{White1993DMRG} computes ground states in MPS form by sweeps of local optimizations, and time-evolving block decimation (TEBD)~\cite{vidal2004tebd} evolves an MPS through local updates followed by truncation. MPS methods have also been applied to classical stochastic dynamics~\cite{johnson2010stochastic,temme2010stochastic}. In numerical analysis, dynamical low-rank approximation projects the evolution onto the tangent space of a manifold of low-rank matrices or tensors~\cite{KochLubich2007,LubichRohwedderSchneiderVandereycken2013}. The projector-splitting integrator carries this out in TT form~\cite{lubich2015tt}. The basis-update and Galerkin integrator was developed for low-rank matrices and extended to Tucker tensors~\cite{CerutiLubich2022}. Low-rank tensor solvers have been developed for chemical master equations in the quantized TT format~\cite{kazeev2014cme}, for HJB equations~\cite{dolgov2021hjb,oster2022hjb}, and for high-dimensional PDEs by combining Monte Carlo with sketching~\cite{chen2025montecarlo}. All of these methods require the equation and solve it within the tensor format, whereas we learn the map from data. The local updates to TT cores in these solvers motivate the layers of TEN-NO.

\paragraph{Tensor networks in machine learning.}
Tensor decompositions have been used to compress the weights of neural networks, with CP~\cite{Lebedev2015} and Tucker~\cite{Kim2016} decompositions of convolution kernels and TT decompositions of fully connected layers~\cite{novikov2015tensorizing}. Low-rank networks can also be trained directly on the low-rank manifold~\cite{Schotthoefer2022}. Tensor networks have also been fitted directly to data, as classifiers~\cite{stoudenmire2016supervised} and as probabilistic models~\cite{tang2025mps,hur2023sketch,peng2025density}. Compositional tensor trains~\cite{EigelMirandaNouySommer2025} compose residual layers whose increments are functional TTs and interpret this composition as a time-discrete dynamical system. Compositions of TTs also underlie transport maps for high-dimensional densities~\cite{CuiDolgov2022} and coordinate flows that reduce TT ranks~\cite{DektorVenturi2023}. In these works, the TT represents the model, that is, the weights of a network or a single function or density. In TEN-NO, a neural network learns a map between input and output states, both represented in TT form.

\subsection{Organization}\label{sec:organization}
Section~\ref{sec:preliminaries} introduces the TT notation and operations. Section~\ref{sec:methods} defines the two TEN-NO variants and analyzes their cost. Section~\ref{sec:experiments} presents the numerical results, and Section~\ref{sec:discussion} discusses the results and their implications for scaling.

\section{Preliminaries}\label{sec:preliminaries}
This section introduces the TT notation and the two TT operations used by the network.

\subsection{TT representation}\label{sec:tt}
For a positive integer $m$, let $[m]=\{1,\ldots,m\}$. We identify a tensor with a function on its index set,
\begin{equation}
  X:[n_1]\times\cdots\times[n_d]\longrightarrow\R,
  \qquad X\in\mathcal V:=\R^{n_1\times\cdots\times n_d}.
  \label{eq:state-space}
\end{equation}
The mode sizes $n_1,\ldots,n_d$ are arbitrary and remain fixed within each learning task. For a spin system, $n_i=2$ and $d$ is the number of spins. For a field on a tensor-product grid, $n_i$ is the number of grid points along the $i$th coordinate and $d$ is the spatial dimension. We write $n=\max_i n_i$ and $x_{a:b}=(x_a,\ldots,x_b)$ for $a\le b$. We use MATLAB notation for slices of multidimensional arrays. When a tensor is viewed as a matrix, the row indices appear before the semicolon and the column indices after it. We use $\langle\cdot,\cdot\rangle_F$ and $\|\cdot\|_F$ for the Frobenius inner product and norm of tensors.

\begin{definition}[TT representation]\label{def:tt}
A tensor $X\in\mathcal V$ has a TT representation with bond dimensions $(r_0,\ldots,r_d)$, where $r_0=r_d=1$, if there exist cores
\begin{equation}
  G_i\in\R^{r_{i-1}\times n_i\times r_i},
  \qquad 1\le i\le d,
  \label{eq:tt-core-shapes}
\end{equation}
such that
\begin{equation}
\begin{aligned}
  X(x_1,\ldots,x_d)
  &=\sum_{\alpha_1=1}^{r_1}\cdots
    \sum_{\alpha_{d-1}=1}^{r_{d-1}}
    \prod_{i=1}^{d}G_i(\alpha_{i-1},x_i,\alpha_i),
  \qquad \alpha_0=\alpha_d=1,\\
  &=G_1(:,x_1,:)\,G_2(:,x_2,:)\cdots G_d(:,x_d,:).
\end{aligned}
\label{eq:tt}
\end{equation}
In the second line, each slice $G_i(:,x_i,:)$ is an $r_{i-1}\times r_i$ matrix. Here $x_i$ is a physical index and $\alpha_i$ is an internal, or bond, index. The smallest bond dimensions for which such cores exist are the TT ranks of $X$.
\end{definition}

\begin{figure}[htbp]
\centering
\begin{tikzpicture}[x=1cm,y=1cm,
  core/.style={circle,draw=black,fill=white,minimum size=5mm,inner sep=0pt},
  bond/.style={line width=.75pt},
  every node/.style={font=\small}]
  \node at (-.95,0) {$\Tt{G}=$};
  \foreach \pos/\name/\idx in {0/g1/1,1.6/g2/2,3.2/g3/3,6.3/gd/d}{
    \node[core] (\name) at (\pos,0) {};
    \draw[bond] (\name.north)--(\pos,.85) node[above] {$x_{\idx}$};
    \node[below=3pt] at (\name.south) {$G_{\idx}$};
  }
  \draw[bond] (g1)--node[above] {$\alpha_1$}(g2);
  \draw[bond] (g2)--node[above] {$\alpha_2$}(g3);
  \draw[bond] (g3)--(4.05,0);
  \node at (4.55,0) {$\cdots$};
  \draw[bond] (5.05,0)--node[above] {$\alpha_{d-1}$}(gd);
  \node[align=left,anchor=west,font=\footnotesize] at (7.15,.3)
    {$x_i\in[n_i]$\\[3pt]$\alpha_i\in[r_i]$};
\end{tikzpicture}
\caption{TT representation of a tensor. The open legs carry the physical indices, and the horizontal bonds are contracted.}
\label{fig:tt-notation}
\end{figure}
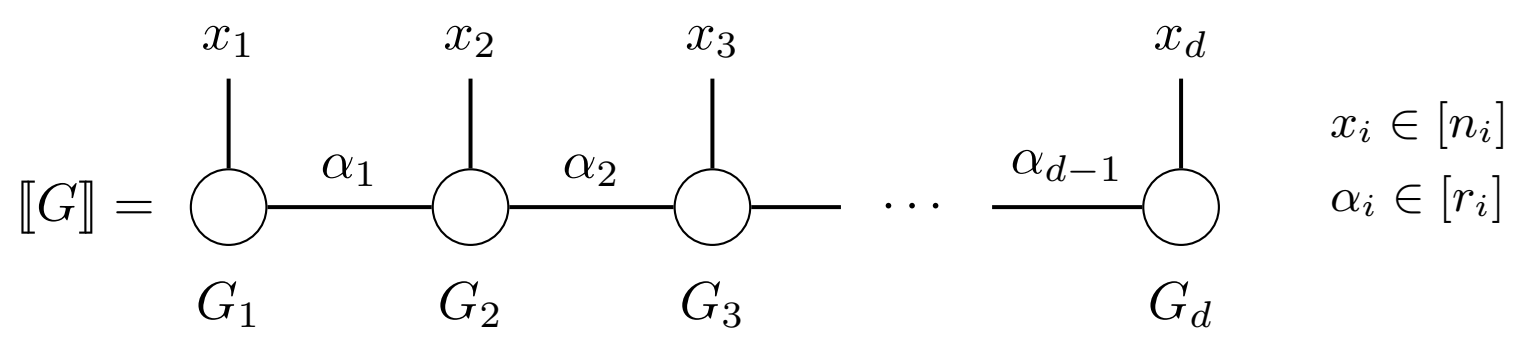

For the ordered core tuple $G=(G_1,\ldots,G_d)$, define the TT with these cores by
\begin{equation}
  \Tt{G}(x_1,\ldots,x_d)
  \coloneqq G_1(:,x_1,:)\,G_2(:,x_2,:)\cdots G_d(:,x_d,:).
  \label{eq:tt-brackets}
\end{equation}
We use the tensor diagrams in Figure~\ref{fig:tt-notation}, which are standard in TT methods. Here, a node denotes a tensor, an open leg denotes a free index, and a line joining two nodes denotes summation over their common index. The trivial bonds $\alpha_0$ and $\alpha_d$ are omitted.

The number of entries in the TT cores satisfies
\begin{equation}
  \sum_{i=1}^{d}n_i r_{i-1}r_i\le dnr^2,
  \label{eq:state-storage}
\end{equation}
where $r$ is any upper bound on the bond dimensions $r_i$. The full tensor has $\prod_{i=1}^{d}n_i$ entries. Each unfolding satisfies $\operatorname{rank}X(x_{1:k};x_{k+1:d})\le r_k$.

\begin{assumption}[Low-rank states]\label{ass:lowrank}
The input and target states admit accurate TT approximations with bond dimensions at most $r_{\mathrm{data}}$.
\end{assumption}

The unfolding rank bound and Oseledets's result~\cite[Theorem~2.1]{oseledets2011} together show that the TT ranks of a tensor equal the ranks of its unfoldings. Thus a tensor has a low-rank TT representation exactly when its unfoldings have low rank. For the problem classes discussed in Section~\ref{sec:introduction}, prior work has reported accurate low-rank TT approximations of probability distributions in stochastic lattice models~\cite{johnson2010stochastic}, ground states of one-dimensional Hamiltonians~\cite{schollwock2011}, and value functions of high-dimensional HJB equations~\cite{dolgov2021hjb}.

\subsection{TT operations}\label{sec:tt-operations}

The network uses two operations that act directly on the TT cores, without forming the full tensor. The Frobenius inner product evaluates the training loss (Section~\ref{sec:loss}), and TT rounding controls the bond dimensions in Algorithm~\ref{alg:ttsvd} and brings the inputs into canonical form (Section~\ref{sec:setup}).

\paragraph{Inner product.}
Let $\Tt{G}$ and $\Tt{H}$ have the same mode sizes. At bond $i$, the index and its dimension are $\alpha_i$ and $r_i$ in $G$, and $\beta_i$ and $r'_i$ in $H$. At each site $i$, contracting $G_i$ and $H_i$ over their shared physical index gives a matrix $M_i$ with paired bond indices. Its entries are given by
\begin{equation}
  M_i(\alpha_{i-1}\beta_{i-1};\alpha_i\beta_i)
  =\sum_{x_i=1}^{n_i}G_i(\alpha_{i-1},x_i,\alpha_i)
                       H_i(\beta_{i-1},x_i,\beta_i).
  \label{eq:tt-inner-core}
\end{equation}
The Frobenius inner product is the contraction of these matrices along their paired bonds~\cite[Section~4.2]{oseledets2011},
\begin{equation}
  \langle \Tt{G},\Tt{H}\rangle_F
  =\sum_{x_1,\ldots,x_d}\Tt{G}(x_1,\ldots,x_d)\,\Tt{H}(x_1,\ldots,x_d)
  =M_1M_2\cdots M_d.
  \label{eq:tt-inner}
\end{equation}
Both ends have bond dimension one, so the product is a scalar. To evaluate the product efficiently, we apply each $M_i$ in factored form. At site $i$, we contract the accumulated result from the preceding sites with the two core slices for each $x_i$ before summing over $x_i$. For bond dimensions bounded by $r$ and $r'$, this costs $O(dn\,rr'(r+r'))$ operations. The Frobenius norm is $\|\Tt{G}\|_F=\langle \Tt{G},\Tt{G}\rangle_F^{1/2}$.

\paragraph{TT-SVD and TT rounding.}
TT-SVD constructs a TT from a full tensor through successive reshapes and truncated SVDs~\cite[Algorithm~1]{oseledets2011}. TT rounding reduces the bond dimensions of an existing TT. A right-to-left QR sweep first puts the cores in right-orthogonal form, and a left-to-right SVD sweep then truncates the bonds~\cite[Algorithm~2]{oseledets2011}. We write $\TTSVD_r$ for rounding with output bond dimensions at most $r$.

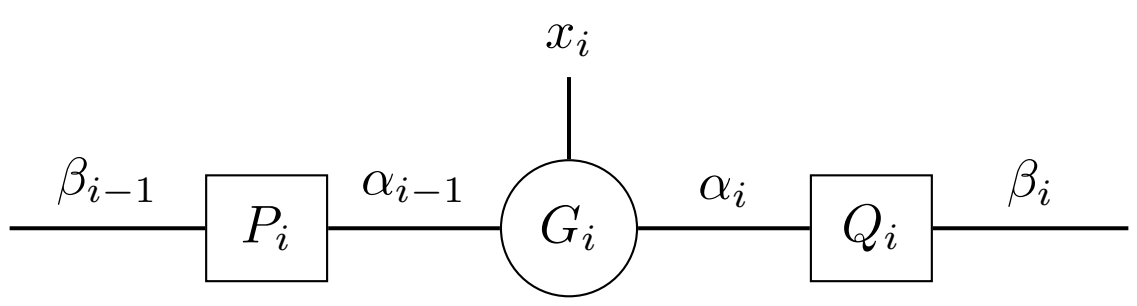
\begin{figure}[H]
\centering
\begin{tikzpicture}[x=1cm,y=1cm,
  core/.style={circle,draw=black,minimum size=9mm,inner sep=0pt},
  factor/.style={draw=black,minimum width=8mm,minimum height=7mm,inner sep=0pt},
  bond/.style={line width=.75pt},
  every node/.style={font=\small}]
  \node[factor] (left) at (-2,0) {$P_i$};
  \node[core] (current) at (0,0) {$G_i$};
  \node[factor] (right) at (2,0) {$Q_i$};
  \draw[bond] (-3.7,0)--node[above] {$\beta_{i-1}$}(left.west)
    (left.east)--node[above] {$\alpha_{i-1}$}(current.west)
    (current.east)--node[above] {$\alpha_i$}(right.west)
    (right.east)--node[above] {$\beta_i$}(3.7,0);
  \draw[bond] (current.north)--(0,1) node[above] {$x_i$};
\end{tikzpicture}
\caption{Bond matrices in the rounded core~\eqref{eq:tt-round-bond-factors}. The open horizontal legs $\beta_{i-1}$ and $\beta_i$ carry the retained bond indices.}
\label{fig:tt-round-bond-factors}
\end{figure}

Let $G=(G_1,\ldots,G_d)$ be the input TT. At site $i$, we apply a matrix $P_i$ to the left bond $\alpha_{i-1}$ of $G_i$ and a matrix $Q_i$ to the right bond $\alpha_i$. The rounded core in Figure~\ref{fig:tt-round-bond-factors} is
\begin{equation}
 \widehat G_i(:,x_i,:)=P_iG_i(:,x_i,:)Q_i,
 \qquad P_1=Q_d=[1].
 \label{eq:tt-round-bond-factors}
\end{equation}

The matrices $P_i$ and $Q_i$ combine the effects of the two sweeps. The sweeps pass factors in opposite directions, so the interior bond matrices can depend on all input cores. For a fixed rank bound and fixed QR and SVD conventions, write
\[
 P_i=\mathcal P_i\bigl(G_1,\ldots,G_d\bigr),
 \qquad
 Q_i=\mathcal Q_i\bigl(G_1,\ldots,G_d\bigr).
\]
The rounded tensor is
\begin{equation}
 \Tt{\widehat G}(x_1,\ldots,x_d)
 =\bigl(P_1G_1(:,x_1,:)Q_1\bigr)
  \cdots
  \bigl(P_dG_d(:,x_d,:)Q_d\bigr).
 \label{eq:tt-round-bond-product}
\end{equation}
The standard SVD sweep returns a left-canonical representation of this same tensor. Appendix~\ref{app:tt-rounding} derives the product form and gives the rounding cost and error bound.

\section{Methods}\label{sec:methods}
Section~\ref{sec:setup} writes the input and target states in TT form. Sections~\ref{sec:local-updates} and~\ref{sec:direct} define the two variants of TEN-NO, and Sections~\ref{sec:loss} and~\ref{sec:cost} give the training objective and the computational cost.

\subsection{Problem setup in TT form}\label{sec:setup}
A fixed preprocessing map $\mathcal C$ represents each state in TT form. It uses TT-SVD for a full tensor and TT rounding for a state already in TT form. The resulting cores are left-canonical, i.e., the unfoldings $G_i(\alpha_{i-1}x_i;\alpha_i)$, $i<d$, have orthonormal columns (Appendix~\ref{app:tt-rounding}), so that the same state always gives the same input cores. For a prescribed bound $r$, we fix the bond dimensions
\begin{equation}
  r_i=\min\Bigl(r,\ \prod_{j\le i}n_j,\ \prod_{j>i}n_j\Bigr),\qquad 0\le i\le d,
  \label{eq:fixed-bond-dims}
\end{equation}
which are the largest TT ranks compatible with the bound $r$. In particular, $r_0=r_d=1$. Both the preprocessing map $\mathcal C$ and the rounding step $\TTSVD_r$ in Algorithm~\ref{alg:ttsvd} produce TTs with exactly these bond dimensions. Each truncation retains $r_i$ singular values. If the rank at bond $i$ is smaller than $r_i$, we complete the left singular vectors to $r_i$ orthonormal columns, assigning zero singular values to the added columns. This completion preserves the represented tensor, keeps the cores left-canonical, and gives them the sizes expected by the local networks. Using the bracket notation in~\eqref{eq:tt-brackets}, the input and target TTs are
\begin{equation}
  \Tt{G^{(0)}}=\mathcal C(\rho_0),
  \qquad \Tt{Y}=\mathcal C(\rho_{\mathrm{tar}}),
  \label{eq:input-cores}
\end{equation}
where $\rho_0$ is the reference state for a parameter-to-state map. At site $i$, the network also receives a feature vector $\omega_i$ derived from $\omega$. The feature vector $\omega_i$ may include global parameters such as the diffusion coefficient. For a fixed-coefficient model, $\omega_i$ is empty.

Let $G^{(\ell)}=(G_i^{(\ell)})_{i=1}^d$ denote the tuple of TT cores after $\ell$ layers. The layer updates and network output are given by
\begin{equation}
  G^{(\ell+1)}=\Phi_\ell\bigl(G^{(\ell)},\omega;\Theta_\ell\bigr),\quad 0\le\ell<L,
  \qquad
  \NN_\Theta(\rho_0,\omega)=\Tt{G^{(L)}},
  \label{eq:network}
\end{equation}
where $\Theta_\ell=(\theta_i^{(\ell)})_{i=1}^d$ collects the local network parameters in layer $\ell$ and $\Theta=(\Theta_0,\ldots,\Theta_{L-1})$. In TEN-NO1, $\Phi_\ell$ applies~\eqref{eq:core-update} at each site and then rounds the resulting TT. In TEN-NO2, $\Phi_\ell$ returns the core tuple $(G_i^{(\ell)}+R_i^{(\ell)})_{i=1}^d$ from the direct update in Section~\ref{sec:direct}. No layer forms a full tensor.

\subsection{Locally connected core update}\label{sec:local-updates}
Let $\mathcal G=(\{1,\ldots,d\},\mathcal E)$ be an undirected graph on the sites. The neighborhood $\mathcal N(i)$ consists of site $i$ and its neighbors in $\mathcal G$,
\begin{equation}
  \mathcal N(i)=\{i\}\cup\bigl\{j:\{i,j\}\in\mathcal E\bigr\}.
  \label{eq:neighborhood}
\end{equation}
The graph describes which sites are coupled and need not follow the order of the cores. For a path with edges $\{i,i+1\}$, an interior site has neighborhood $\mathcal N(i)=\{i-1,i,i+1\}$. Adding the edge $\{1,d\}$ gives periodic neighborhoods. In a two-dimensional lattice, neighboring sites can be far apart in the TT ordering. At site $i$ in layer $\ell$, the local network $\NN^{(\ell)}_i$ takes the neighborhood cores $G^{(\ell)}_{\mathcal N(i)}=(G^{(\ell)}_j)_{j\in\mathcal N(i)}$ and the corresponding site feature vectors $\omega_{\mathcal N(i)}=(\omega_j)_{j\in\mathcal N(i)}$ as input. With trainable parameters $\theta^{(\ell)}_i$, it returns a residual core $R_i^{(\ell)}$
\begin{equation}
  R^{(\ell)}_i=\NN^{(\ell)}_i\bigl(G^{(\ell)}_{\mathcal N(i)},\omega_{\mathcal N(i)};\theta^{(\ell)}_i\bigr)\in\R^{\tilde r_{i-1}\times n_i\times\tilde r_i},
  \qquad \tilde r_i\ge r_i,\quad \tilde r_0=\tilde r_d=1.
  \label{eq:residual}
\end{equation}
The bond dimension $r_i$ of the incoming cores $G^{(\ell)}$ is given by~\eqref{eq:fixed-bond-dims}, and $\tilde r_i$ is the corresponding bond dimension before rounding. The input and output sizes of each local network are fixed and shared across samples. All residuals of layer $\ell$ are computed from the cores $G^{(\ell)}$. The updated core $\widetilde G_i^{(\ell)}$ is
\begin{equation}
  \widetilde G^{(\ell)}_i=\iota_i\bigl(G^{(\ell)}_i\bigr)+R^{(\ell)}_i,
  \label{eq:core-update}
\end{equation}
where the zero-padding map $\iota_i$ places $G^{(\ell)}_i$ in the leading $r_{i-1}\times n_i\times r_i$ block of an $\tilde r_{i-1}\times n_i\times\tilde r_i$ array and sets the remaining entries to zero. The padded cores represent the same tensor as $G^{(\ell)}$, and $\iota_i$ is the identity if $\tilde r_{i-1}=r_{i-1}$ and $\tilde r_i=r_i$. Any nonlinear activation is part of $\NN^{(\ell)}_i$, and the incoming core enters~\eqref{eq:core-update} through the identity branch $\iota_i(G^{(\ell)}_i)$. After forming all updated cores $\widetilde G_i^{(\ell)}$, we round the TT they represent using $\TTSVD_r$ (Section~\ref{sec:tt-operations}). The resulting cores $G^{(\ell+1)}$ have the bond dimensions in~\eqref{eq:fixed-bond-dims} required by the local networks in the next layer. Rounding truncates only at bonds with $\tilde r_i>r_i$. Algorithm~\ref{alg:ttsvd} summarizes this variant.

\begin{algorithm}[H]
\hrule\smallskip
\caption{Locally connected network with TT-SVD (TEN-NO1)}\label{alg:ttsvd}
\smallskip\hrule\smallskip
\small
\noindent\textbf{Input:} initial cores $G^{(0)}_1,\ldots,G^{(0)}_d$; neighborhoods $\mathcal N(i)$; site feature vectors $\omega_1,\ldots,\omega_d$; local networks $\NN^{(\ell)}_i$ with trainable parameters $\theta^{(\ell)}_i$, $0\le\ell<L$, $1\le i\le d$; residual bond dimensions $\tilde r_i$; bond-dimension bound $r$.\\
\textbf{Output:} cores $G^{(L)}_1,\ldots,G^{(L)}_d$.
\begin{enumerate}[label=\arabic*:,leftmargin=2.2em,itemsep=1pt,topsep=3pt]
\item \textbf{for} $\ell=0,\ldots,L-1$ \textbf{do}
\item \hspace*{1.2em}\textbf{for} $i=1,\ldots,d$ \textbf{do}
\item \hspace*{2.4em}$R^{(\ell)}_i\gets\NN^{(\ell)}_i\bigl(G^{(\ell)}_{\mathcal N(i)},\omega_{\mathcal N(i)};\theta^{(\ell)}_i\bigr)$
\item \hspace*{2.4em}$\widetilde G^{(\ell)}_i\gets\iota_i\bigl(G^{(\ell)}_i\bigr)+R^{(\ell)}_i$
\item \hspace*{1.2em}\textbf{end for}
\item \hspace*{1.2em}$\bigl(G^{(\ell+1)}_1,\ldots,G^{(\ell+1)}_d\bigr)\gets\TTSVD_r\bigl(\widetilde G^{(\ell)}_1,\ldots,\widetilde G^{(\ell)}_d\bigr)$
\item \textbf{end for}
\item \textbf{return} $\bigl(G^{(L)}_1,\ldots,G^{(L)}_d\bigr)$
\end{enumerate}
\hrule
\end{algorithm}

TEBD applies two-site gates to an MPS and truncates the affected bond after each gate~\cite{vidal2004tebd,unfried2023tebd}. This pattern of local updates followed by rank control motivates Algorithm~\ref{alg:ttsvd}, in which the local maps are learned from data and TT rounding acts on the updated tensor at the end of each layer.

The residual formulation lets each layer learn a correction to the current state~\cite{he2016resnet}. When all residual cores $R_i^{(\ell)}$ are zero, the padded cores represent the input tensor. Since its TT rank at each bond $i$ is at most $r_i$, rounding preserves that tensor.

\subsection{Direct core update}\label{sec:direct}
In TEN-NO1, the input $G$ to TT rounding in Section~\ref{sec:tt-operations} is the chain $\widetilde G^{(\ell)}$ formed by the residual update in~\eqref{eq:core-update}. By~\eqref{eq:tt-round-bond-factors}, the rounded TT has a representation with core slices $P_i\widetilde G_i^{(\ell)}(:,x_i,:)Q_i$. The bond matrices $P_i$ and $Q_i$ can depend on all $d$ provisional cores.

For TEN-NO2, consider the bond maps $\mathcal P_i$ and $\mathcal Q_i$ on the incoming TT $G^{(\ell)}$. For the states considered here, we assume that the resulting bond matrices can be approximated from the incoming neighborhood cores $G_{\mathcal N(i)}^{(\ell)}$,
\begin{equation}
\begin{aligned}
 P_i &= \mathcal P_i\bigl(G_1^{(\ell)},\ldots,G_d^{(\ell)}\bigr)
       \approx \widehat P_i\bigl(G_{\mathcal N(i)}^{(\ell)}\bigr),\\
 Q_i &= \mathcal Q_i\bigl(G_1^{(\ell)},\ldots,G_d^{(\ell)}\bigr)
       \approx \widehat Q_i\bigl(G_{\mathcal N(i)}^{(\ell)}\bigr).
\end{aligned}
 \label{eq:local-bond-factors}
\end{equation}
TEN-NO2 uses this assumption to model the effect of the residual update and the two bond multiplications in a single local update. The MLP in~\eqref{eq:residual} reads the incoming neighborhood cores and site features $\omega_{\mathcal N(i)}$ and predicts a residual $R_i^{(\ell)}$ of size $r_{i-1}\times n_i\times r_i$. The updated core $G_i^{(\ell+1)}=G_i^{(\ell)}+R_i^{(\ell)}$ therefore has the prescribed bond dimensions and passes directly to the next layer without TT rounding. Algorithm~\ref{alg:direct} updates the sites in parallel because every MLP reads from the incoming layer.

\begin{algorithm}[H]
\hrule\smallskip
\caption{Direct locally connected network (TEN-NO2)}\label{alg:direct}
\smallskip\hrule\smallskip
\small
\noindent\textbf{Input:} initial cores $G^{(0)}_1,\ldots,G^{(0)}_d$; neighborhoods $\mathcal N(i)$; site feature vectors $\omega_1,\ldots,\omega_d$; local networks $\NN^{(\ell)}_i$ with trainable parameters $\theta^{(\ell)}_i$, $0\le\ell<L$, $1\le i\le d$, whose outputs have the same size as $G^{(\ell)}_i$, that is, $\tilde r_i=r_i$.\\
\textbf{Output:} cores $G^{(L)}_1,\ldots,G^{(L)}_d$.
\begin{enumerate}[label=\arabic*:,leftmargin=2.2em,itemsep=1pt,topsep=3pt]
\item \textbf{for} $\ell=0,\ldots,L-1$ \textbf{do}
\item \hspace*{1.2em}\textbf{for} $i=1,\ldots,d$ \textbf{do}
\item \hspace*{2.4em}$R^{(\ell)}_i\gets\NN^{(\ell)}_i\bigl(G^{(\ell)}_{\mathcal N(i)},\omega_{\mathcal N(i)};\theta^{(\ell)}_i\bigr)$
\item \hspace*{2.4em}$G^{(\ell+1)}_i\gets G^{(\ell)}_i+R^{(\ell)}_i$
\item \hspace*{1.2em}\textbf{end for}
\item \textbf{end for}
\item \textbf{return} $\bigl(G^{(L)}_1,\ldots,G^{(L)}_d\bigr)$
\end{enumerate}
\hrule
\end{algorithm}

The interactions in our examples motivate the choice of neighborhood in~\eqref{eq:local-bond-factors}. In open-chain Glauber dynamics, the flip rate $w_i$ depends only on the spin at site $i$ and the spins at its immediate neighbors. The initial distribution $\rho_0$ is a product distribution (Section~\ref{sec:glauberexperiment}). In the HJB example, the drift component $v_i$ couples only coordinate $i$ and its adjacent coordinates (Section~\ref{sec:hjbexperiment}).

Neighborhood tensor updates for projected entangled pair states (PEPS) use a similar approximation of the environment, evaluating bond truncation errors on a small cluster of tensors around the updated bond~\cite{dziarmaga2021neighborhood,sinha2024pepsmetts,sokolov2025truncating}.

\subsection{Training objective}\label{sec:loss}
For sample $m$, write $G^{(0),m}$ and $Y^{(m)}$ for the input and reference cores from Section~\ref{sec:setup}. Given the input cores $G^{(0),m}$ and the problem parameters $\omega^{(m)}$, the network returns the predicted cores $G^{(L),m}$. For evolution maps, we minimize $\mathcal J(\Theta)$, the mean squared relative error over the $N$ training samples,
\begin{equation}
  \mathcal J(\Theta)=\frac1N\sum_{m=1}^{N}
  \frac{\bigl\|\Tt{G^{(L),m}}-\Tt{Y^{(m)}}\bigr\|_F^2}{\bigl\|\Tt{Y^{(m)}}\bigr\|_F^2}.
  \label{eq:loss}
\end{equation}
We evaluate the loss by TT contractions. With $G=G^{(L),m}$ and $H=Y^{(m)}$, the numerator of the $m$th term in~\eqref{eq:loss} expands as
\begin{equation}
  \bigl\|\Tt{G}-\Tt{H}\bigr\|_F^2=\langle \Tt{G},\Tt{G}\rangle_F-2\langle \Tt{G},\Tt{H}\rangle_F+\langle \Tt{H},\Tt{H}\rangle_F.
  \label{eq:loss-inner}
\end{equation}
Each inner product in~\eqref{eq:loss-inner} is evaluated using the contraction in~\eqref{eq:tt-inner}. The last term in~\eqref{eq:loss-inner} is also the denominator $\|\Tt{Y^{(m)}}\|_F^2$ in~\eqref{eq:loss} for sample $m$. We assume $\|\Tt{Y^{(m)}}\|_F>0$ so that the relative error is well defined.

\subsection{Computational cost}\label{sec:cost}
For $j=1,2$, let $C_{\NN,j}$ be an upper bound on the cost of evaluating one local network in variant $j$. The bond dimensions are bounded by $r$ for the network states, by $\tilde r$ for the states before rounding, and by $r_{\mathrm{data}}$ for the target TTs. Table~\ref{tab:costs} compares the two variants.

\begin{table}[htbp]
\centering
\caption{Arithmetic costs of an $L$-layer forward pass and loss evaluation, together with the number of TT core entries saved for backpropagation. All quantities are reported per sample.}
\label{tab:costs}
\begin{tabular*}{\linewidth}{@{\extracolsep{\fill}}lll@{}}
\toprule
Step & TEN-NO1 (Algorithm~\ref{alg:ttsvd}) & TEN-NO2 (Algorithm~\ref{alg:direct})\\
\midrule
Local networks & $O(LdC_{\NN,1})$ & $O(LdC_{\NN,2})$\\
Core updates & $O(Ldn\tilde r^2)$ & $O(Ldnr^2)$\\
Rounding & $O(Ldn\tilde r^3)$ & ---\\
Loss & $O(dn(r+r_{\mathrm{data}})^3)$ & $O(dn(r+r_{\mathrm{data}})^3)$\\
Saved TT cores for backpropagation & $O(Ldn\tilde r^2)$ & $O(Ldnr^2)$\\
\bottomrule
\end{tabular*}
\end{table}

The bounds in Table~\ref{tab:costs} grow at most linearly with the number of sites $d$ for fixed network depth $L$, maximum mode size $n$, and bond-dimension bounds, provided that neighborhood sizes, feature dimensions, and local network widths are bounded. In TEN-NO1, the rounding step in each layer contributes the cubic term in $\tilde r$ and requires a sequence of $2(d-1)$ dependent factorizations. In TEN-NO2, the $d$ site updates within a layer are independent. For local MLPs with one hidden layer of fixed width $w$, $C_{\NN,2}=O(wnr^2)$ when the neighborhood sizes and feature dimensions are bounded. Under these conditions, only the loss evaluation in the TEN-NO2 column of Table~\ref{tab:costs} has a cubic arithmetic cost bound in the bond dimensions. Appendix~\ref{app:cost} gives the derivation.

\section{Numerical results}\label{sec:experiments}
We test fixed-time evolution maps for the HJB equation and open-chain Glauber dynamics, and a parameter-to-state map for the periodic TFI model.
Both TEN-NO variants use $L=3$ layers with a separate MLP $\NN_i^{(\ell)}$ at each site and layer. Each MLP takes the concatenated TT core entries and site features $\omega_j$ as input and has one hidden layer of width 64 with SiLU activation.

We use 2,048 training, 256 validation, and 256 regular test samples at each reported problem size. We train in float32 on an NVIDIA GeForce RTX 4060 Laptop GPU using AdamW and a OneCycleLR schedule. In the scaling studies, we use an effective batch size of 128 and allow up to 10,000 updates. Each larger TEN-NO2 model is warm-started from the preceding checkpoint with a fresh optimizer and schedule. For each reported size, the cumulative training time $T_{\mathrm{C}}$ is the sum of the wall-clock times spent training and validating that model and all preceding models in its warm-start chain. Training stops when the validation error plateaus or the update budget is reached. We evaluate the checkpoint with the lowest validation error.

\subsection{HJB equation with local drift}\label{sec:hjbexperiment}
The prediction target is the terminal field $u(T,\cdot)$ at $T=1$. On the periodic domain $[-\pi,\pi)^d$, the equation is
\begin{equation}
\partial_t u=\nu\Delta u+v(x)\cdot\nabla u-\norm{\nabla u}_2^2,
\qquad u(0,x)=u_0(x).
\label{eq:hjb}
\end{equation}
Here $\nu$ is the diffusion coefficient. The local MLPs use path-graph neighborhoods $\mathcal N(i)$ as defined in Section~\ref{sec:local-updates}. The drift component $v_i$ couples coordinate $i$ to the other sites in $\mathcal N(i)$ with strength $c_i$,
\begin{equation}
v_i(x)=0.2\sin x_i\left(1+c_i
\prod_{j\in\mathcal N(i)\setminus\{i\}}\cos x_j\right).
\end{equation}
We use the following family of initial fields
\[
u_0(x)=0.1\sum_{i=1}^d a_i\cos(x_i+\phi_i).
\]
The amplitudes and phases are mutually independent, with $a_i\sim\operatorname{Uniform}[0.5,1.5]$ and $\phi_i\sim\operatorname{Uniform}[-\pi,\pi]$. The grid has $n_i=32$ points per coordinate.

For TEN-NO, we represent each sampled initial field as a TT and zero-pad its cores to the network bond dimensions with $r=4$. Reference solutions computed with a PyTorch implementation of TT-based RK4 serve as targets for the loss in~\eqref{eq:loss}, and we report the mean relative Frobenius error on the test set. In each of the fixed- and varying-coefficient experiments, we assessed the reliability of the reference solutions by recomputing the terminal fields for four samples at each tested dimension, using a finer grid and smaller time steps. The relative Frobenius differences were below $5\times10^{-4}$, supporting the reliability of the reference solutions used for evaluation.

\subsubsection{Time evolution with fixed coefficients}\label{sec:hjb-fixed}
The coefficients are $c_i=0.1$ for each coordinate and $\nu=0.05$. At $d=3$, we compare TEN-NO1 and TEN-NO2 with FNO~\cite{li2021fno}. For the two TEN-NO variants, the input TTs are put in left-canonical form (Section~\ref{sec:setup}), and their core norms are then balanced for numerical stability. All three runs use the same training, validation, and test pairs on a $32^3$ grid, a batch size of 32, and a budget of 3,000 updates. The two TEN-NO variants start from identical MLP weights and use the same bond dimensions throughout, with 543,552 trainable parameters each. TEN-NO1 uses automatic differentiation through its per-layer TT-SVD. FNO takes the grid values of $u_0$ as input and has 658,897 trainable parameters.

\begin{table}[H]
\centering
\caption{Three-dimensional fixed-coefficient HJB benchmark. Both TEN-NO variants use canonical inputs.}
\label{tab:hjb-three-dimensional}
\begin{tabular*}{\linewidth}{@{\extracolsep{\fill}}lrr@{}}
\toprule
Model & Test error (\%) & Training time (s)\\
\midrule
TEN-NO1 & 0.9548 & 148.14\\
TEN-NO2 & 0.7291 & 25.34\\
FNO & 0.4769 & 448.59\\
\bottomrule
\end{tabular*}
\end{table}

In the three-dimensional benchmark, TEN-NO2 yields a lower test error and requires less training time than TEN-NO1 (Table~\ref{tab:hjb-three-dimensional}). TEN-NO2 achieves a test error comparable to FNO's with a $17.7\times$ training speedup. In TT form, the inputs and reference solutions require a factor of $23.3$ less storage than the dense arrays used by FNO. At this resolution, a dense-grid FNO requires $32^d$ values per feature channel, which limits its use in higher dimensions. We therefore use TEN-NO2 in the remaining experiments, starting with the fixed-coefficient HJB problem at $d=5,10,15,$ and $20$ (Table~\ref{tab:continuousresults}(a)).

\begin{table}[H]
\centering
\caption{(a) Fixed-coefficient HJB predictions with $c_i=0.1$ and $\nu=0.05$.}
\label{tab:continuousresults}
\begin{tabular*}{\linewidth}{@{\extracolsep{\fill}}lrr@{}}
\toprule
Dimension $d$ & Test error (\%) & $T_{\mathrm{C}}$ (s)\\
\midrule
5 & 1.1829 & 32.37\\
10 & 1.4423 & 86.27\\
15 & 1.2731 & 118.60\\
20 & 1.2671 & 160.84\\
\bottomrule
\end{tabular*}
\end{table}

The mean test error remains around $1.3\%$ across the tested dimensions (Table~\ref{tab:continuousresults}(a)). For $d=15$ and $d=20$, we lower the peak learning rate from $3\times10^{-3}$ to $10^{-3}$. Both runs meet the validation plateau criterion shortly after the learning rate peaks, so each adds less training time to the chain than the $d=10$ run.

\subsubsection{Time evolution with varying coefficients}
The problem parameters are $\omega=(\boldsymbol c,\nu)$. The site features $\omega_i$ encode the local coupling $c_i$ and the shared diffusion coefficient $\nu$. We set $c_i=sq_i$, with $q_i$ obtained by rescaling independent $\operatorname{Uniform}[0.75,1.25]$ draws to have mean one across sites. The training grid spans $s\in[0.05,0.20]$ and $\nu\in[0.035,0.065]$.

The test set combines $(s,\nu)$ pairs from the training grid with held-out initial fields and normalized coupling vectors $\boldsymbol q$. For interpolation and extrapolation, we reuse these fields and vectors and vary one scalar at a time, evaluating each new value at both endpoints of the other scalar's training range. Interpolation uses midpoints between training grid values, and extrapolation uses $s\in\{0,0.25\}$ or $\nu\in\{0.03,0.07\}$. The interpolation and extrapolation columns in Table~\ref{tab:continuousparameterized}(b) pool the errors from varying $s$ and $\nu$.

\begin{table}[!htbp]
\ContinuedFloat
\centering
\caption{(b) HJB predictions with heterogeneous couplings and a varying diffusion coefficient.}
\label{tab:continuousparameterized}
\begin{tabular*}{\linewidth}{@{\extracolsep{\fill}}lrrrr@{}}
\toprule
Dimension $d$ & Test error (\%) & Interpolation (\%) & Extrapolation (\%) & $T_{\mathrm{C}}$ (s)\\
\midrule
5 & 0.5721 & 0.5839 & 0.7234 & 28.53\\
10 & 0.3569 & 0.3649 & 0.4698 & 97.11\\
15 & 0.4199 & 0.4243 & 0.5227 & 176.70\\
20 & 0.4468 & 0.4539 & 1.5996 & 278.16\\
\bottomrule
\end{tabular*}
\end{table}

Test and interpolation errors remain close, reaching a minimum at $d=10$ and rising modestly thereafter. Extrapolation error follows a similar trend through $d=15$ and shows its largest increase at $d=20$. The training time added at each warm-start stage increases with dimension. TEN-NO2 scales well with dimension in the varying-coefficient setting.

\subsection{Glauber dynamics on an open chain}\label{sec:glauberexperiment}
We learn the fixed-time map $\rho_0\mapsto\rho_T$ for the open-chain kinetic Ising model~\cite{glauber1963} at $T=0.25$. A spin configuration $s\in\{-1,+1\}^d$ has energy
\begin{equation}
E(s)=-\sum_{i=1}^{d-1}J_i s_i s_{i+1}-\sum_{i=1}^{d}b_i s_i.
\end{equation}
At inverse temperature $\beta$, the heat-bath rate for flipping spin $i$ is
\begin{equation}
w_i(s)=\frac12\left[1-s_i\tanh\!\left(\beta\left[J_{i-1}s_{i-1}+J_i s_{i+1}+b_i\right]\right)\right],
\end{equation}
where terms involving $s_0$ or $s_{d+1}$ are omitted. Let $s^{(i)}$ denote the configuration obtained by flipping spin $i$. The distribution $\rho_t$ satisfies
\begin{equation}
\partial_t\rho_t(s)=\sum_{i=1}^{d}
\left[w_i(s^{(i)})\rho_t(s^{(i)})-w_i(s)\rho_t(s)\right].
\label{eq:glauber}
\end{equation}
The local MLP at site $i$ uses the same path-graph neighborhood $\mathcal N(i)$ as the flip rate $w_i$ (Section~\ref{sec:local-updates}).

The initial distributions have the product form
\begin{equation}
\rho_0(s)=\prod_{i=1}^{d} g_i(s_i),
\qquad g_i(s_i)=\frac{1+a_i s_i}{2}.
\label{eq:glauber-product-input}
\end{equation}
The $a_i$ are drawn independently from $\operatorname{Uniform}[-0.8,0.8]$. Each product distribution is represented as a rank-one TT, and its cores are zero-padded to dimension $r=5$ on each internal bond.

Reference distributions $\rho_T$ computed in TT form with a PyTorch implementation of TEBD~\cite{vidal2004tebd} serve as targets for the loss in~\eqref{eq:loss}, and we report the mean relative Frobenius error on the test set. To assess the reliability of the reference distributions, we recomputed three samples from the fixed-coefficient experiment and five from the varying-coefficient experiment at each chain length, using half the time step and twice the TT rank bound. The largest relative Frobenius difference was about $1.21\times10^{-5}$, supporting the reliability of the reference distributions used for evaluation.

\subsubsection{Time evolution with fixed coefficients}
With $\beta=0.8$, $J_i=1$, and $b_i=0.5$ fixed, we test TEN-NO2 on progressively longer chains (Table~\ref{tab:glauberresults}(a)).

\begin{table}[!htbp]
\centering
\caption{(a) Fixed-coefficient Glauber predictions at $T=0.25$.}
\label{tab:glauberresults}
\begin{tabular*}{\linewidth}{@{\extracolsep{\fill}}lrr@{}}
\toprule
Chain length $d$ & Test error (\%) & $T_{\mathrm{C}}$ (s)\\
\midrule
32 & 0.7898 & 93.48\\
64 & 0.7555 & 178.10\\
128 & 1.0959 & 338.17\\
256 & 1.6941 & 657.23\\
\bottomrule
\end{tabular*}
\end{table}

Cumulative training time grows approximately linearly with chain length (Table~\ref{tab:glauberresults}(a)). The test error changes little between $d=32$ and $d=64$ and increases at the two larger sizes.

\subsubsection{Time evolution with varying coefficients}
The inverse temperature is fixed at $\beta=1$. For problem parameters $\omega=(\boldsymbol J,\boldsymbol b)$, the site features $\omega_i$ encode the local field $b_i$ and the adjacent couplings $J_{i-1}$ and $J_i$. We draw the base couplings $J_i^0$ independently from $\operatorname{Uniform}[1,2]$ and set $J_i=\lambda J_i^0$ and $b_i=\lambda h_0$. The training grid spans $\lambda\in[0.1,0.4]$ and $h_0\in[0.1,0.5]$.

The test set uses held-out initial distributions and base couplings, with $(\lambda,h_0)$ taken from the training grid. Interpolation and extrapolation use a common subset of the test initial distributions and a common subset of the training base couplings. The $h_0$ values are also shared, so only $\lambda$ varies. Extrapolation uses $\lambda\in\{0.05,0.45\}$. Table~\ref{tab:glauberparameterized}(b) reports the errors and cumulative training time by chain length.

\begin{table}[H]
\ContinuedFloat
\centering
\caption{(b) Glauber predictions with heterogeneous couplings.}
\label{tab:glauberparameterized}
\begin{tabular*}{\linewidth}{@{\extracolsep{\fill}}lrrrr@{}}
\toprule
Chain length $d$ & Test error (\%) & Interpolation (\%) & Extrapolation (\%) & $T_{\mathrm{C}}$ (s)\\
\midrule
32 & 0.4071 & 0.3513 & 0.7182 & 95.24\\
64 & 0.6188 & 0.5485 & 1.1534 & 178.56\\
128 & 0.9452 & 0.8192 & 1.7496 & 343.87\\
256 & 1.4375 & 1.2407 & 2.8919 & 669.57\\
\bottomrule
\end{tabular*}
\end{table}

All three errors increase with chain length, and on the matched tests, extrapolation error exceeds interpolation error by a growing margin. Cumulative training time is approximately proportional to chain length (Table~\ref{tab:glauberparameterized}(b)).

\subsection{Ground states of the TFI model}\label{sec:tfiexperiment}
Given the transverse-field vector $\omega=\boldsymbol{h}=(h_1,\ldots,h_d)$, we predict the ground-state MPS $\rho_{\mathrm{tar}}=\psi_{\mathrm{gs}}(\boldsymbol{h})$ of a periodic TFI chain. The Hamiltonian is
\begin{equation}
H(\boldsymbol{h})=-\sum_{i=1}^d Z_iZ_{i+1}-\sum_{i=1}^d h_i X_i,
\qquad Z_{d+1}=Z_1,
\end{equation}
where $X_i$ and $Z_i$ are Pauli operators at site $i$. We draw the $h_i$ independently from $\operatorname{Uniform}[0.7,1.3]$.

At each chain length, all samples share a precomputed ground-state MPS $\rho_0$ for the uniform field $\boldsymbol{h}=\boldsymbol{1}$, rounded to the network's TT rank bound $r=20$. Using the site features $\omega_i$ derived from the sampled fields, the network predicts corrections to the cores of the shared MPS $\rho_0$. At site $i$, the local MLP reads field features from the periodic neighborhood $\mathcal N(i)$ and cores from the contiguous window of that neighborhood in the open TT (Section~\ref{sec:direct}). For each sampled $\boldsymbol{h}$, we compute the reference ground-state MPS with the DMRG solver in \href{https://github.com/ITensor/ITensorMPS.jl/tree/v0.3.20}{\texttt{ITensorMPS.jl} v0.3.20}~\cite{White1993DMRG,FishmanWhiteStoudenmire2022}. We assessed the reliability of each reference MPS by checking even parity and bounding its phase-aligned distance to the ground state using the eigenstate residual and the even-sector gap. All numerical bounds were below $10^{-3}$, supporting the reliability of the reference states.

\subsubsection{Ground-state prediction across chain lengths}

For a given $\boldsymbol{h}$, let $\psi=\psi_{\mathrm{gs}}(\boldsymbol{h})$ denote the reference state and $\widehat\psi$ denote the prediction. The normalized overlap is
\[
\eta=\frac{\bigl\lvert \langle \psi,\widehat{\psi}\rangle \bigr\rvert}
          {\bigl\lVert\psi\bigr\rVert_2\,\bigl\lVert\widehat{\psi}\bigr\rVert_2}.
\]
We compute $\eta$ by MPS contraction and use $2-2\eta$ as the per-sample training loss. The test error is the mean phase-aligned Frobenius error between unit-normalized states, with per-sample error $\sqrt{2-2\eta}$. We evaluate TEN-NO2 on periodic chains of 16, 32, and 48 spins (Table~\ref{tab:tfiresults}).

\begin{table}[H]
\centering
\caption{Ground-state predictions for the periodic TFI chain using a TT rank bound of $20$.}
\label{tab:tfiresults}
\begin{tabular*}{\linewidth}{@{\extracolsep{\fill}}lrr@{}}
\toprule
Chain length $d$ & Test error (\%) & $T_{\mathrm{C}}$ (s)\\
\midrule
16 & 0.6368 & 90.94\\
32 & 1.6237 & 314.12\\
48 & 3.1529 & 1111.95\\
\bottomrule
\end{tabular*}
\end{table}

The test error increases with chain length, and the cumulative training time grows faster than linearly in the number of spins, especially from $d=32$ to $d=48$ (Table~\ref{tab:tfiresults}). In additional TT-SVD tests on the reference MPSs, the TT rank needed for comparable rounding accuracy grew approximately linearly with chain length. This may partly explain the rise in TEN-NO2's test error at a fixed rank bound.

\section{Discussion}\label{sec:discussion}
We introduce TEN-NO to learn high-dimensional solution maps from an input state and problem parameters. A dense state with $d$ modes of size $n$ requires $n^d$ entries. Keeping states in low-rank TT form throughout the network makes prediction feasible at sizes where dense storage is impractical.

At each site, a local network predicts a residual update from the TT cores and site features in its graph neighborhood. TEBD's use of truncation after local updates motivates TEN-NO1, which rounds the updated TT at the end of each layer. TEN-NO2 assumes that the bond matrices $P_i$ and $Q_i$ can be approximated from the incoming neighborhood cores, as in~\eqref{eq:local-bond-factors}. This motivates direct updates at fixed bond dimensions, with the site updates able to run in parallel.

In the three-dimensional HJB benchmark, TEN-NO2 reaches a test error comparable to FNO's in substantially less training time. Its layers process compact TT cores, reducing the amount of state data handled relative to FNO's dense-grid features. This may help explain the timing difference. Under the conditions in Section~\ref{sec:cost}, the forward pass and loss evaluation each cost $O(d)$ arithmetic operations per sample. Prediction accuracy depends on whether the states admit accurate low-rank TT approximations and whether the selected neighborhoods provide enough information for the core updates. The TFI rounding tests showed that longer chains required higher ranks for comparable rounding accuracy. This motivates testing adaptive ranks at a prescribed accuracy and measuring how the resulting cost scales.

\appendix
\section{TT rounding}\label{app:tt-rounding}
Let $G=(G_1,\ldots,G_d)$ be the input TT, with mode sizes at most $n$ and bond dimensions at most $\tilde r$. At a truncation step $i$, let $L_i$ and $S_i$ be the matrices formed by contracting the current cores strictly to the left and right of site $i$, respectively. Their rows index the corresponding physical multi-indices, and their columns index the adjacent bonds. With $A_i$ denoting the left unfolding of the current core and $x_i$ varying fastest in the row indices of both unfoldings, the current tensor $X$ satisfies
\[
 X(x_{1:i};x_{i+1:d})=(L_i\otimes I_{n_i})A_i S_i^{\mathsf T}.
\]
When $L_i$ and $S_i$ have orthonormal columns and $A_i=U\Sigma V^{\mathsf T}$ is a compact SVD, the full unfolding has the compact SVD $[(L_i\otimes I_{n_i})U]\Sigma(S_iV)^{\mathsf T}$.

\paragraph{Right-to-left orthogonalization.}
Set $C_d=[1]$. For $i=d,\ldots,2$, unfold the slices $G_i(:,x_i,:)C_i$ into a matrix $B_i$ with the left bond as its row index. Take the thin QR factorization $B_i^{\mathsf T}=W_iR_i$, reshape $W_i^{\mathsf T}$ into the swept core $H_i$, and set $C_{i-1}=R_i^{\mathsf T}$. At the first site set $H_1(:,x_1,:)=G_1(:,x_1,:)C_1$ and $C_0=[1]$. Thus
\[
 G_i(:,x_i,:)C_i=C_{i-1}H_i(:,x_i,:),\qquad 1\le i\le d.
\]
The sweep preserves the represented tensor and gives right-orthogonal cores $H_2,\ldots,H_d$~\cite[Lemma~3.1]{oseledets2011}.

\paragraph{Left-to-right truncation.}
For $i=1,\ldots,d-1$, compute the truncated SVD of the current core unfolding $A_i$, keeping the leading $r$ singular values or all available values if there are fewer,
\begin{equation}
 A_i\approx U_i\Sigma_iV_i^{\mathsf T}.
 \label{eq:tt-round-svd}
\end{equation}
Reshape $U_i$ to replace the current core and multiply $\Sigma_iV_i^{\mathsf T}$ into the left bond of the next core. Before step $i$, $L_i$ has orthonormal columns because the preceding cores have been replaced by orthonormal factors. The cores to the right of site $i$ retain their form from the right-to-left sweep, so $S_i$ also has orthonormal columns. The factorization above then shows that this step truncates the SVD of the full unfolding of the current tensor. After the sweep, the first $d-1$ left unfoldings have orthonormal columns. The factorizations in both sweeps act on matrices with at most $n\tilde r$ rows and $\tilde r$ columns, so each costs $O(n\tilde r^3)$ operations and the two sweeps cost $O(dn\tilde r^3)$.

\paragraph{Bond matrices.}
Start with $P_1=[1]$. If the factor entering site $i$ in the SVD sweep is $P_iC_{i-1}$, the QR relation gives $A_i=D_iC_i$, where $D_i$ is the left unfolding of the slices $P_iG_i(:,x_i,:)$. For $i<d$, set
\[
 P_{i+1}=U_i^{\mathsf T}D_i,\qquad Q_i=C_iV_i\Sigma_i^\dagger,
\]
where $\Sigma_i^\dagger$ is the diagonal pseudoinverse. Then $P_{i+1}C_i=U_i^{\mathsf T}A_i=\Sigma_iV_i^{\mathsf T}$, which is the factor passed to the next core. Also, $D_iQ_i=A_iV_i\Sigma_i^\dagger=U_i\Sigma_i\Sigma_i^\dagger$. A retained zero singular value gives a zero row in the transfer and hence a zero slice on the incoming bond of the next product core. Zeroing the corresponding columns of $U_i$ therefore leaves the contracted tensor unchanged. At the last site, the final core slices are $P_dG_d(:,x_d,:)$ because $C_d=[1]$. Taking $Q_d=[1]$ gives the product~\eqref{eq:tt-round-bond-product} with the core slices in~\eqref{eq:tt-round-bond-factors}.

\paragraph{Error.}
Let $X^{(0)}$ be the input tensor and $X^{(i)}$ the tensor after step $i$ of the truncation sweep, and let $E_i=X^{(i-1)}-X^{(i)}$. The columns of the $i$th unfolding of $E_i$ lie in the range of $(L_i\otimes I_{n_i})(I-U_iU_i^{\mathsf T})$, which is orthogonal to the range of the new left interface $L_{i+1}=(L_i\otimes I_{n_i})U_i$. Later steps leave the first $i$ truncated cores fixed, so for $j>i$ the columns of the $i$th unfolding of $E_j$ lie in the range of $L_{i+1}$. Hence $\langle E_i,E_j\rangle_F=0$ for $i\ne j$, and
\[
 \bigl\|X^{(0)}-X^{(d-1)}\bigr\|_F^2=\sum_{i=1}^{d-1}\|E_i\|_F^2,
\]
where $\|E_i\|_F^2$ is the sum of the squared singular values discarded in~\eqref{eq:tt-round-svd}. Let $\varepsilon_i$ denote the minimum Frobenius error when approximating the $i$th unfolding of $X^{(0)}$ by a matrix of rank at most $r$. The current tensor $X^{(i-1)}$ is obtained from $X^{(0)}$ by an orthogonal projection in the variables $x_{1:i-1}$. Since orthogonal projections do not increase singular values, $\|E_i\|_F\le\varepsilon_i$, and hence $\|X^{(0)}-X^{(d-1)}\|_F^2\le\sum_{i=1}^{d-1}\varepsilon_i^2$~\cite[Theorem~2.2]{oseledets2011}.

\section{Computational cost}\label{app:cost}
Table~\ref{tab:op-costs} lists the cost of the operations used by the network.

\begin{table}[htbp]
\centering
\caption{Cost of the operations used by the network.}
\label{tab:op-costs}
\begin{tabular*}{\linewidth}{@{\extracolsep{\fill}}ll@{}}
\toprule
Operation & Cost or storage\\
\midrule
Storage of a TT with bond dimensions at most $r$ & $O(dnr^2)$ numbers\\
One local network $\NN^{(\ell)}_i$ in variant $j$, \eqref{eq:residual} & $C_{\NN,j}$\\
One core update~\eqref{eq:core-update} & $O(n\tilde r^2)$\\
TT rounding of a TT with bond dimensions at most $\tilde r$ & $O(dn\tilde r^3)$\\
Inner product of TTs with bond dimensions $r$ and $r'$, \eqref{eq:tt-inner} & $O(dn\,rr'(r+r'))$\\
\bottomrule
\end{tabular*}
\end{table}

\paragraph{Algorithm~\ref{alg:ttsvd}.}
A layer performs two steps.
\begin{enumerate}[label=(\alph*)]
\item \emph{Local updates.} The $d$ local networks cost $O(dC_{\NN,1})$ operations, and forming the cores~\eqref{eq:core-update} costs $O(dn\tilde r^2)$.
\item \emph{Rounding.} The QR sweep and the truncated SVD sweep act on cores of size at most $\tilde r\times n\times\tilde r$, at a cost of $O(n\tilde r^3)$ operations per core, that is, $O(dn\tilde r^3)$ in total (Appendix~\ref{app:tt-rounding}). Each factorization uses the output of the previous one, so this step is a chain of $2(d-1)$ dependent factorizations.
\end{enumerate}
The per-layer cost is therefore
\begin{equation}
  C_{1}=O\bigl(dC_{\NN,1}+dn\tilde r^3\bigr).
  \label{eq:cost-alg1}
\end{equation}

\paragraph{Algorithm~\ref{alg:direct}.}
A layer performs only step (a), with $\tilde r=r$, so its cost is
\begin{equation}
  C_{2}=O\bigl(dC_{\NN,2}+dnr^2\bigr).
  \label{eq:cost-alg2}
\end{equation}
Compared with~\eqref{eq:cost-alg1}, the cubic term $dn\tilde r^3$ of the rounding is absent. The $d$ site updates within a layer can run in parallel, so the number of sequential steps is independent of $d$. In Algorithm~\ref{alg:ttsvd}, the rounding step requires a sequence of factorizations whose length grows linearly with $d$.

\paragraph{Local networks.}
Suppose each local network is an MLP with one hidden layer of width $w$. Let $s=\max_i|\mathcal N(i)|$, and let $f$ bound the number of scalar features in $\omega_{\mathcal N(i)}$. Each local network receives at most $snr^2+f$ input entries. By~\eqref{eq:residual}, its output contains at most $n\tilde r^2$ entries in TEN-NO1 and $nr^2$ in TEN-NO2. The two dense maps and the hidden activation therefore give
\[
 C_{\NN,1}=O\bigl(w(snr^2+f+n\tilde r^2)\bigr),
 \qquad
 C_{\NN,2}=O\bigl(w(snr^2+f+nr^2)\bigr).
\]
For bounded $s$ and $f$, the bound for TEN-NO2 simplifies to $O(wnr^2)$.

\paragraph{Loss.}
Computing the three inner products in~\eqref{eq:loss-inner} costs $O(dn(r+r_{\mathrm{data}})^3)$ operations per sample (Table~\ref{tab:op-costs}). The last inner product also gives the denominator $\|\Tt{Y^{(m)}}\|_F^2$.

\paragraph{Overall cost.}
For fixed $n$, $r$, $\tilde r$, and $r_{\mathrm{data}}$ and bounded $s$, $f$, and $w$, the per-layer bounds in~\eqref{eq:cost-alg1} and~\eqref{eq:cost-alg2} and the loss bound are all $O(d)$. The forward pass and loss together cost $O(Ld)$ operations per sample, which is linear in $d$ for fixed $L$. One TT state requires $O(dnr^2)$ entries. Counting the TT cores saved across all $L$ layers for backpropagation gives $O(Ldnr^2)$ entries in Algorithm~\ref{alg:direct} and $O(Ldn\tilde r^2)$ in Algorithm~\ref{alg:ttsvd}. For fixed $L$ and local network sizes, site-specific networks have a parameter count that grows linearly with $d$, and networks shared across sites have a count independent of $d$. A dense representation of a single state requires $\prod_{i=1}^{d}n_i$ entries.

\end{document}